\documentclass{article}
\usepackage{amssymb, amsthm, enumerate}
\usepackage{graphicx}
\usepackage{indentfirst}

\newtheorem{tetel}{Theorem}

\newtheorem{excercise}{Excercise}

\begin{document}

\title{On convex closed planar curves as equidistant sets}
\markright{On convex closed planar curves as equidistant sets}
\author{Csaba Vincze\\
Institute of Mathematics, University of Debrecen,\\ P.O.Box 400, H-4002 Debrecen, Hungary\\
csvincze@science.unideb.hu}
\footnotetext[2]{{\bf MR subject classification:} 51M04}
\footnotetext[3]{The work is supported by the EFOP-3.6.1-16-2016-00022 project. The project is co-financed by the European Union and the European Social Fund.}
\maketitle

\begin{abstract}
The equidistant set of two nonempty subsets $K$ and $L$ in the Euclidean plane is a set all of whose points have the same distance from $K$ and $L$. Since the classical conics can be also given in this way, equidistant sets can be considered as a kind of their generalizations: $K$ and $L$ are called the focal sets. In their paper \cite{PS} the authors posed the problem of the characterization of closed subsets in the Euclidean plane that can be realized as the equidistant set of two connected disjoint closed sets. We prove that any convex closed planar curve can be given as an equidistant set, i.e. the set of equidistant curves contains the entire class of convex closed planar curves. In this sense the equidistancy is a generalization of the convexity. 
\end{abstract}

\section{Introduction.}

Let $K\subset \mathbb{R}^2$ be a subset in the Euclidean coordinate plane. The distance between a point $X\in \mathbb{R}^2$ and $K$ is measured by the usual infimum formula:
$$d(X, K) := \inf  \{ d(X,A) \ | \ A\in K \},$$
where $d(X,A)$ is the Euclidean distance of the point $X$ and $A$ running through the elements of $K$. Let us define the equidistant set of $K$ and $L\subset \mathbb{R}^2$ as the set all of whose points have the same distance from $K$ and $L$:
$$\{K=L\}:=\{ X\in \mathbb{R}^2\ | \ d(X, K)=d(X, L)\}.$$
The equidistant sets can be considered as a kind of the generalization of conics \cite{PS}: $K$ and $L$ are called the focal sets. Equidistant sets are often called midsets. Their investigations have been started by Wilker's and Loveland's fundamental works \cite{Wilker} and \cite{Loveland}. For another generalization of the classical conics and their applications see
e.g. \cite{ErdosVincze}, \cite{MF} (polyellipses and their applications), \cite{GS}, \cite{VA_1}, \cite{VA_2}, \cite{VA_3} and \cite{VA_4}. Let $R>0$ be a positive real number. The \emph{parallel body} of a set $K\subset \mathbb{R}^2$ with radius $R$ is the union of the closed disks with radius $R$ centered at the points of $K$. The infimum of the positive numbers such that $L$ is a subset of the parallel body of $K$ with radius $R$ and vice versa is called 
the \emph{Hausdorff distance} of $K$ and $L$. It is well-known that the Hausdorff metric makes the family of non\-empty closed and bounded (i.e. compact) subsets in the plane a complete metric space. Our main result is the characterization of convex closed curves in the plane as equidistant sets. The result is a contribution to the open problem posed in \cite{PS}: \emph{characterize all closed sets of the plane that can be realized as the equidistant set of two connected disjoint closed sets.}

\begin{tetel}
Any convex closed curve in the plane can be given as an equidistant set.
\end{tetel}

\section{The proof of Theorem 1.}

The proof is divided into two parts. The first step is to prove the statement for convex polygons\footnote{In general there are lots of different ways to give a convex polygon as an equidistant set. We are going to construct finite focal sets that seems to be the best choice for the computer simulation \cite{VA_4}.}. In the second step we use a continuity argument based on the following theorem. 

\begin{tetel}  {\emph{\cite{PS}}} If $K$ and $L$ are disjoint compact subsets in the plane, $K_n\to K$ and $L_n\to L$ are convergent sequences of nonempty compact subsets with respect to the Hausdorff metric then for any $R>0$ we have that 
$$\{K_n=L_n\}\cap \bar{D}(R) \ \to \ \{K=L\}\cap \bar{D}(R),$$
where $\bar{D}(R)$ denotes the closed disk with radius $R$ centered at the origin
\end{tetel}

\begin{figure}
\centering
\includegraphics[scale=0.43]{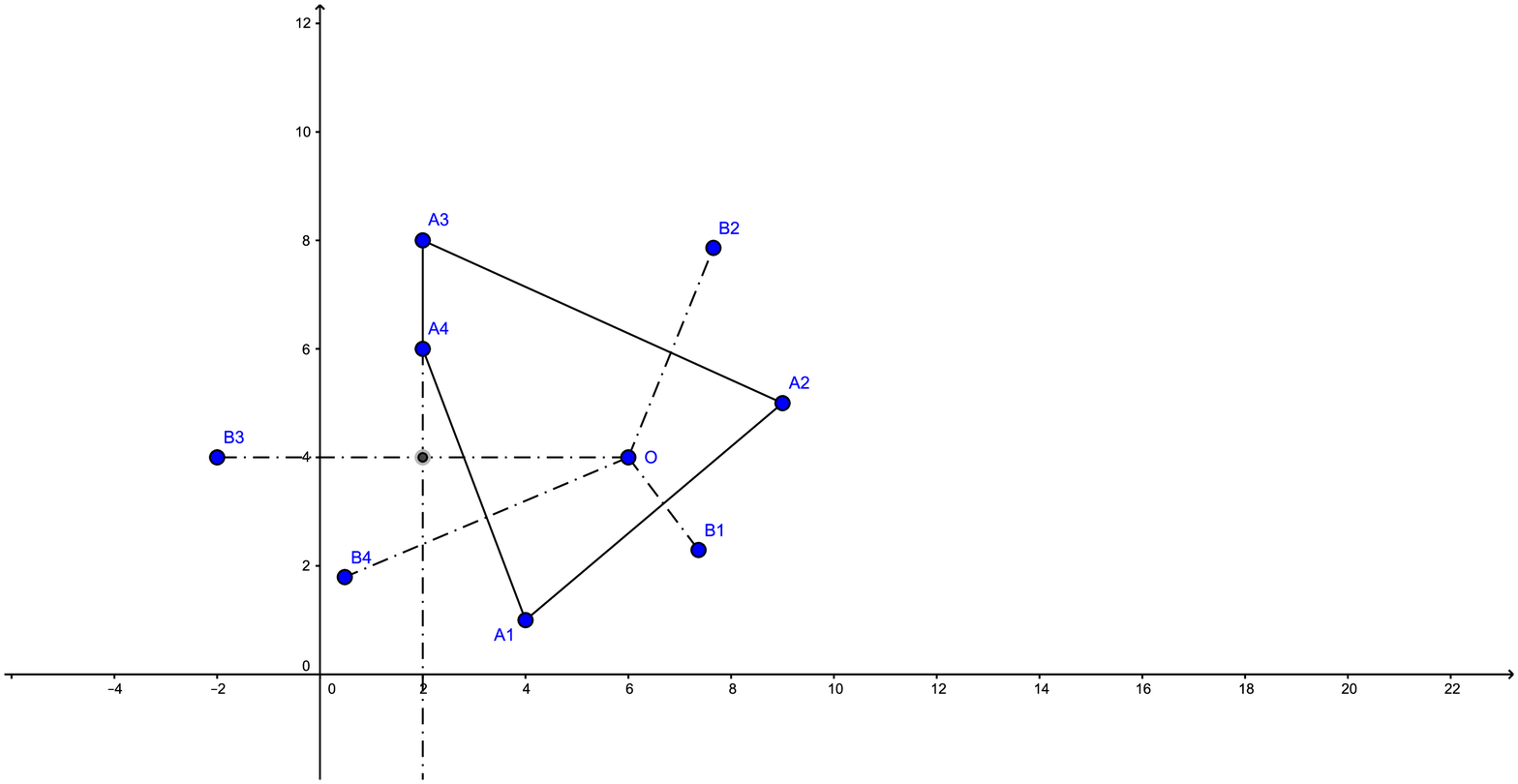}
\caption{The reflected pairs.}
\end{figure}

\subsection{The first step: convex polygons as equidistant sets.} Let $P$ be a convex polygon in the plane with adjacent vertices $A_1$, $\ldots$, $A_n$, $A_{n+1}:=A_1$, $A_{n+2}:=A_2$, $\ldots,$ (in counterclockwise direction), i.e. the indices are taken by modulo $n$. Let us denote the edges by $a_{ii+1}:=A_iA_{i+1}\ \ (i=1, \ldots, n).$
After choosing a point $O$ in the interior of the convex hull of $P$, consider the reflected pairs $(O, B_i)\ \ (i=1, \ldots, n),$ where the supporting line of the edge $a_{ii+1}$ is the perpendicular bisector of $OB_i$ for any $i=1, \ldots, n$ (see figure 1). Let us denote the Voronoi cells with respect to the set of points  $B_1, \ldots, B_n$ by $V_1, \ldots, V_n$. They give the decomposition of the plane into non-empty closed convex subsets such that for any $X\in V_i$ 
$$d(X, B_i)=\min_{1 \leq j \leq n} d(X, B_j)\ \ (i=1, \ldots, n).$$
Using the convexity of the polygon we have that $B_1, \ldots, B_n$ are pairwise different points because $B_1=B_j$ ($j\neq 1$) implies that  the edges $a_{12}$ and $a_{jj+1}$ has a common supporting line $l_{1j}$. Since $l_{1j}$ is the supporting line to the entire convex hull of $P$ at the same time it follows that $a_{12}$ and $a_{jj+1}$ are contained in the same segment $a_{1j+1}$ of the polygon which is obviously a contradiction. It can be easily seen that 
$$V_i=\cap_{j\neq i} F_{ij}\ \ (i=1, \ldots, n),$$
where the half plane $F_{ij}$ is bounded by the perpendicular bisector of the segments $B_iB_j$ such that $B_i\in F_{ij}$ for any $j\neq i$ (see figure 2). First of all note that $A_{i+1}$ can not be in the interior of $V_i$ because it is lying on the boundary of the half plane $F_{ii+1}$ for any $i=1, \ldots, n$. Indeed, $B_{i}$ and $B_{i+1}$ are the reflected pairs of the same point $O$ about the edges $a_{ii+1}$ and $a_{i+1i+2}$. Since they are adjacent at $A_{i+1}$ we have
\begin{equation}
d(A_{i+1}, B_i)=d(A_{i+1},O)=d(A_{i+1}, B_{i+1}).
\end{equation}

\begin{figure}
\centering
\includegraphics[scale=0.43]{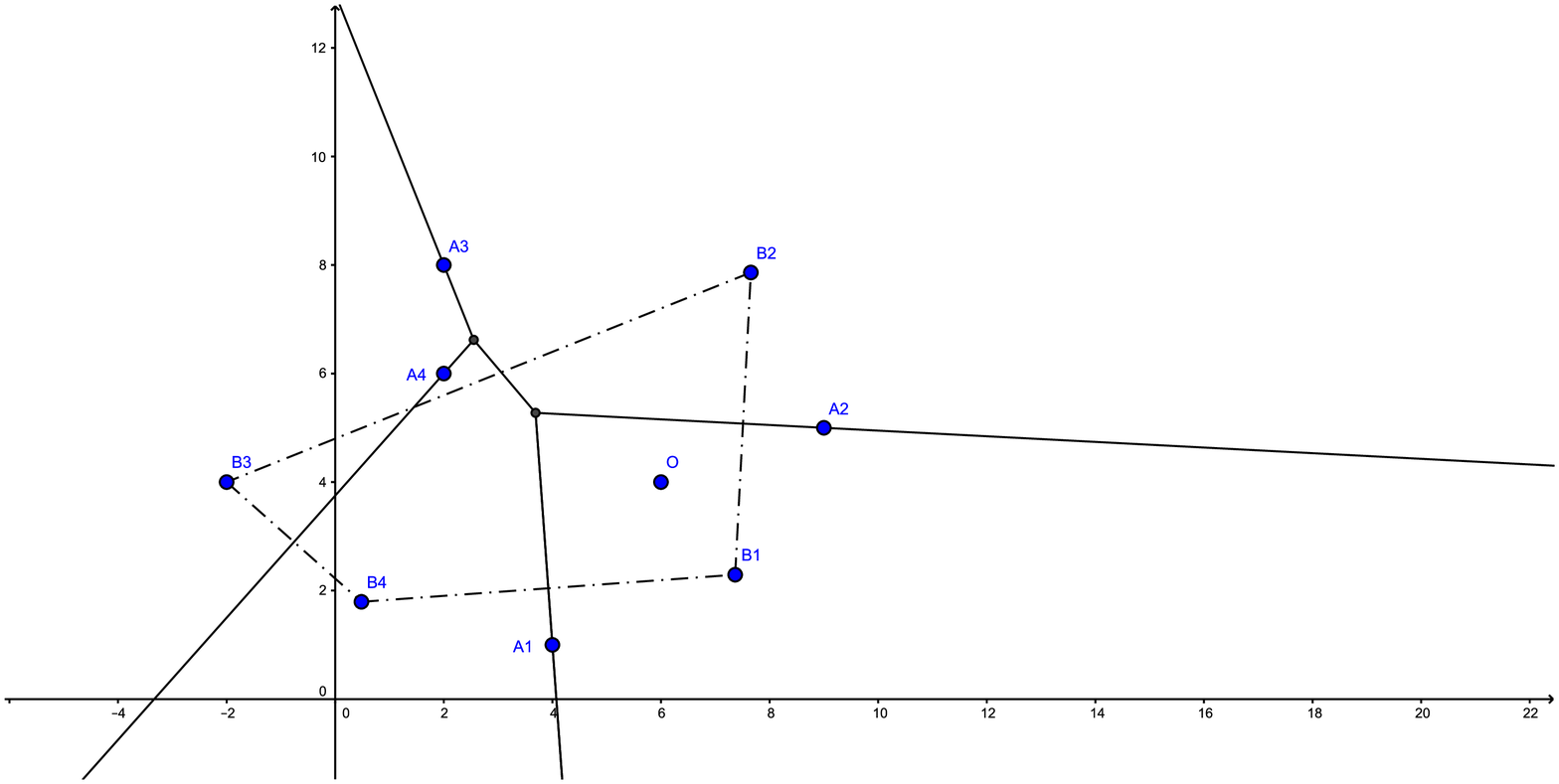}
\caption{The Voronoi cells.}
\end{figure}

Suppose that $A_{i+1}$ is a point outside the Voronoi cell $V_i$. This means that there is a $j\neq i$ such that the perpendicular bisector of $B_i B_{j}$ strictly separates the points $A_{i+1}$ and $B_i$. It is not too hard to see that the points $B_i$, $O$ and $B_j$ can not be collinear in this case because the collinearity $B_i - O - B_j$ implies the edges $A_iA_{i+1}$, $A_jA_{j+1}$ and the perpendicular bisector of $B_iB_j$ to be parallel such that the perpendicular bisector is between the (parallel) supporting lines of the edges. Therefore the strict separation of the points $A_{i+1}$ and $B_i$ is impossible with the perpendicular bisector of $B_iB_j$. Consider now the triangle $B_iOB_j$; the perpendicular bisectors $l_i$, $l_j$ and $l$ belonging to the sides $B_iO$, $B_jO$ and $B_iB_j$ meet the center $C$ of the circumscribed circle (see figure 3). Since $l_i$ and $l_j$ are supporting lines of some edges of the polygon it follows that $P$ is contained in the intersection of the half planes containing $O$ bounded by $l_i$ and $l_j$, respectively. If $l$ strictly separates $A_{i+1}$ (lying on $l_i$) and $B_i$ then $A_{i+1}$ is contained in the open half plane determined by $l_j$ but opposite to the half plane containing $O$. This is a contradiction because $A_{i+1}\in P$. 
We have just proved that $A_{i+1}$ is on the boundary of the Voronoi cell $V_i$. In a similar way we have that $A_i$ is on the boundary of the Voronoi cell $V_i$. Therefore for any $X\in a_{ii+1}$
$$d(X,L)=d(X, B_i)=d(X,O)=d(X,K),$$
where $K:=\{O\}\ \ \textrm{and}\ \ L:=\{B_1, \ldots, B_n\}.$ This means that $P$ is the subset of the equidistant points to $K$ and $L$. Conversely, if $X$ is equidistant to $K$ and $L$ then it is contained in the intersection of a Voronoi cell $V_i$ and the perpendicular bisector of $B_iO$ for some index $i=1, \ldots, n$. Therefore $X$ is lying on the edge $A_iA_{i+1}$.

\begin{figure}
\centering
\includegraphics[scale=0.4]{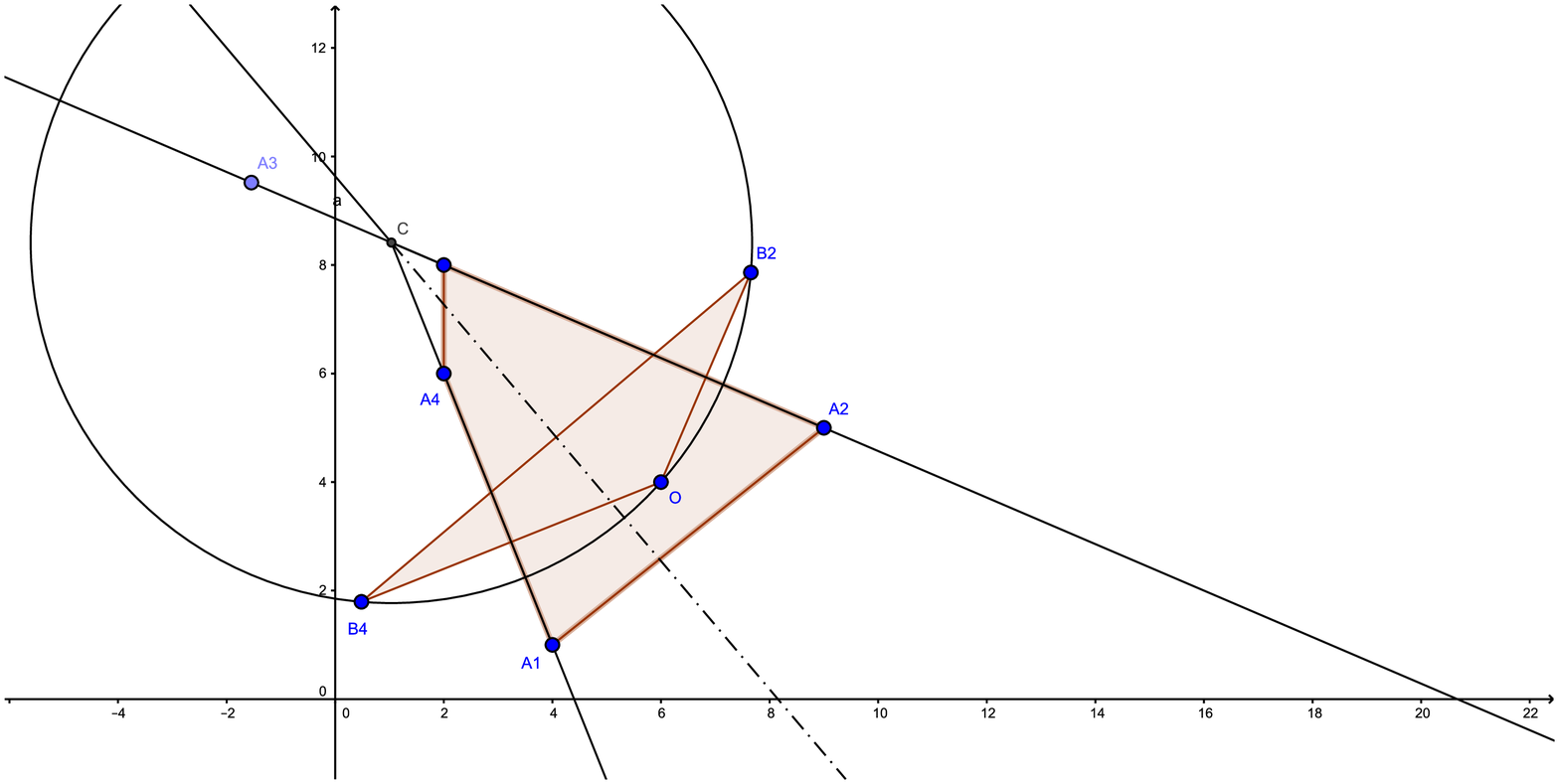}
\caption{The contradiction under the choice $i=2$ and $j=4$.}
\end{figure}

\subsection{The second step: a continuity argument.} Let $C$ be a closed convex curve in the plane and consider its approximation by a sequence of inscribed convex $n$ - gons ($n\geq 3$) with respect to the Hausdorff distance:
$\lim_{n\to \infty}P_n =C.$ As we have seen in the first step $P_n$ can be supposed to be the equidistant set of 
$$K:=\{O\}\ \ \textrm{and}\ \ L_n:=\{B_{1n}, \ldots, B_{nn}\},$$
where $O$ is choosen independently of $n$ in the interior of $C$. It can be easily seen that there exist positive numbers $r$ and $R$ such that $L_n$ is in the closed ring bounded by the circles of radiuses $r$ and $R$ around the point $O$. The collection of compact subsets in such a ring equipped with the Hausdorff metric forms a totally bounded\footnote{ A metric space $(M,d)$ is totally bounded if, for any $\varepsilon > 0$, there exists a finite subset
$F:=\{F_1, \ldots, F_m\}\subset M$ such that $d(X,F):=\min_{1\leq i \leq m} d(X, F_i) < \varepsilon$ for all $X\in M$. In the special case of the metric space of compact subsets in the ring choose a finite collection $\mathcal{N}$ of points such that the union of the open disks centered at the points in $\mathcal{N}$ with radius $\varepsilon$ covers the ring; $F$ can be choosen as the (finite set) of all subsets in $\mathcal{N}$.} complete metric space because of the completeness of the metric space of non-empty compact subsets in $\mathbb{R}^2$  with respect to the Hausdorff metric. It is known that totally boundedness and completeness is equivalent to the sequentially compactness for metric spaces. This means that we can suppose the sequence $L_n$ to be convergent by choosing a convergent subsequence if necessary. By the continuity theorem (Theorem 2) of equidistant sets, the set containing the equidistant points to $K$ and $L$ is the Hausdorff limit of $P_n$ as the sequence of the equidistant sets to $K$ and $L_n$ tending to $L$, i.e. $C$ is the equidistant set of $K$ and $L$. 

\begin{excercise}
Prove that if the convex polygon $P$ with $n$ vertices has an inscribed circle of radius $r$ then it can be given as the equidistant set of
$$K:=\{O\}\ \ \textrm{and}\ \ L:=\{B_1, \ldots, B_n\},$$
where the points $B_1, \ldots, B_n$ are the vertices of a convex polygon having a circumscribed circle of radius $R=2r$ and these circles are concentric.
\end{excercise}
Hint. Choose the center of the inscribed circle of $P$ as the point $O$ in the first step.

\begin{excercise}
What about the higher dimensional version of Theorem 1?
\end{excercise}

\section{Concluding remarks}

\begin{figure}
\centering
\includegraphics[scale=0.3]{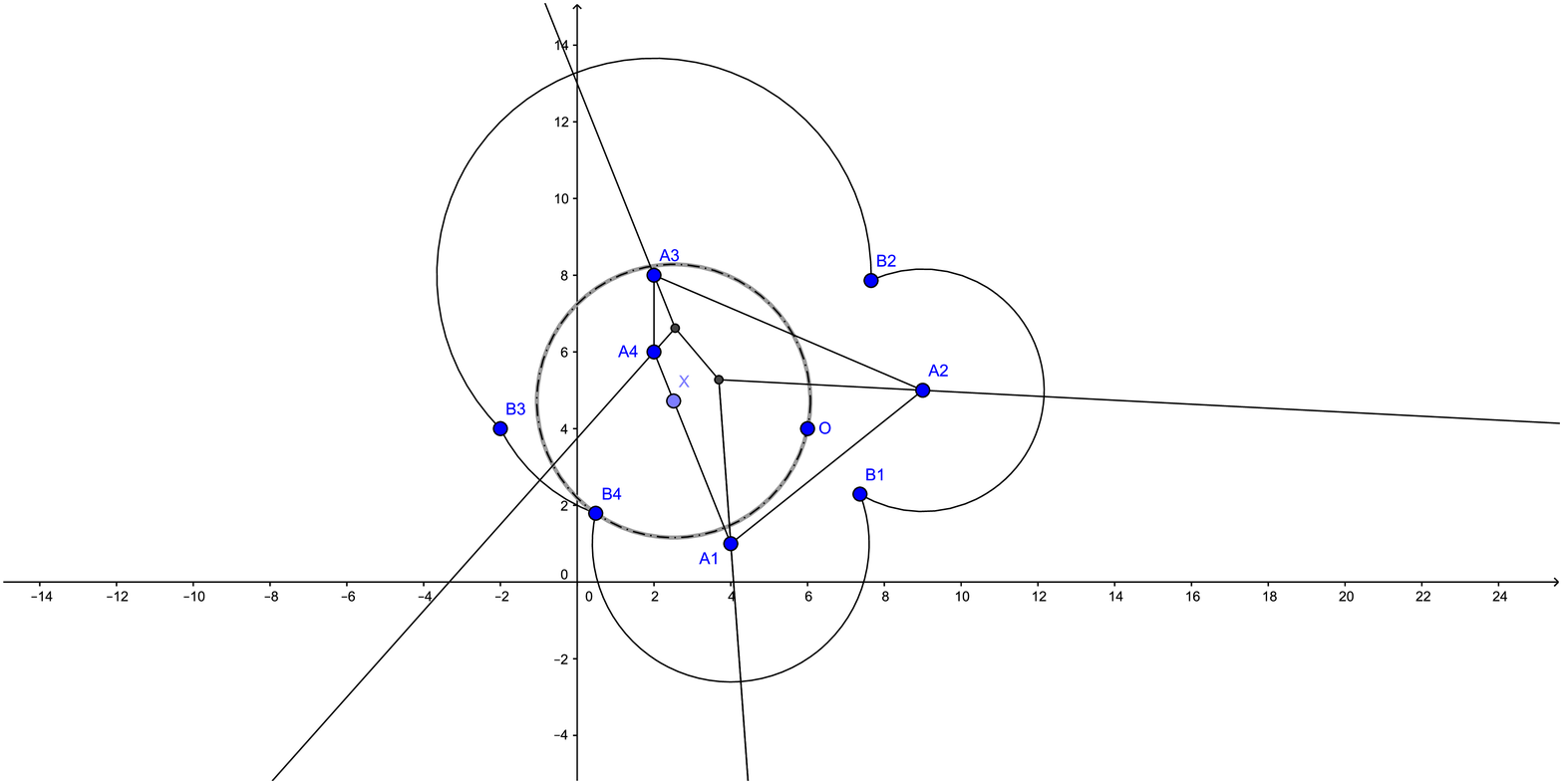}
\caption{The circular arcs in the Voronoi cells.}
\end{figure}

To prove that the class of convex closed planar curves belong to the class of closed subsets in the Euclidean plane that can be realized as the equidistant set of two connected disjoint closed sets we need to present a connected focal set $M$ in the first step of the proof instead of the finite pointset $L$. In each Voronoi cell $V_i$ consider the union $M_i$ of circular arcs passing through $B_{i}$ with centers $A_{i}$ and $A_{i+1}$, respectively (see figure 4).
It can be easily seen that the domain bounded by the union of these "double" arcs (as $i$ runs through $1$ to $n$) contains any circle passing through $B_i$ with center $X\in a_{i i+1}$ as $X$ runs through the points of the polygon $P$; cf. figure $4$ under the choice $i=4$. Therefore if
$$M:=\bigcup_{i=1}^n M_i$$  
then for any $X\in a_{ii+1}$
$$d(X,M)=d(X, B_i)\ \ (i=1, \ldots, n),$$
i.e. the polygon belongs to the equidistant set of $K=\{O\}$ and $M$. Conversely, if $X$ is equidistant to $K$ and $M$
then there exists a point $B\in M$ such that $X$ is lying on the perpendicular bisector of $OB$. For the sake of simplicity consider figure 5, where the positions of the points $X$ and $B$ are illustrated in. Since each circle determined by the arcs passes through the point $O$, it follows that 
$$\angle X A_2 B =\angle X A_2 O < \angle X A_2 B_1.$$
On the other hand $d(A_2,B)=d(A_2,B_1)$ and $A_2X$ is a common side of the triangles $A_2 X B$ and $A_2 X B_1$. Therefore 
$d(X,B_1) < d(X,B),$ i.e. $d(X,M)=d(X,B_1)$ and $X$ is lying on the edge $a_{12}$ of the polygon.

\begin{figure}
\centering
\includegraphics[scale=0.3]{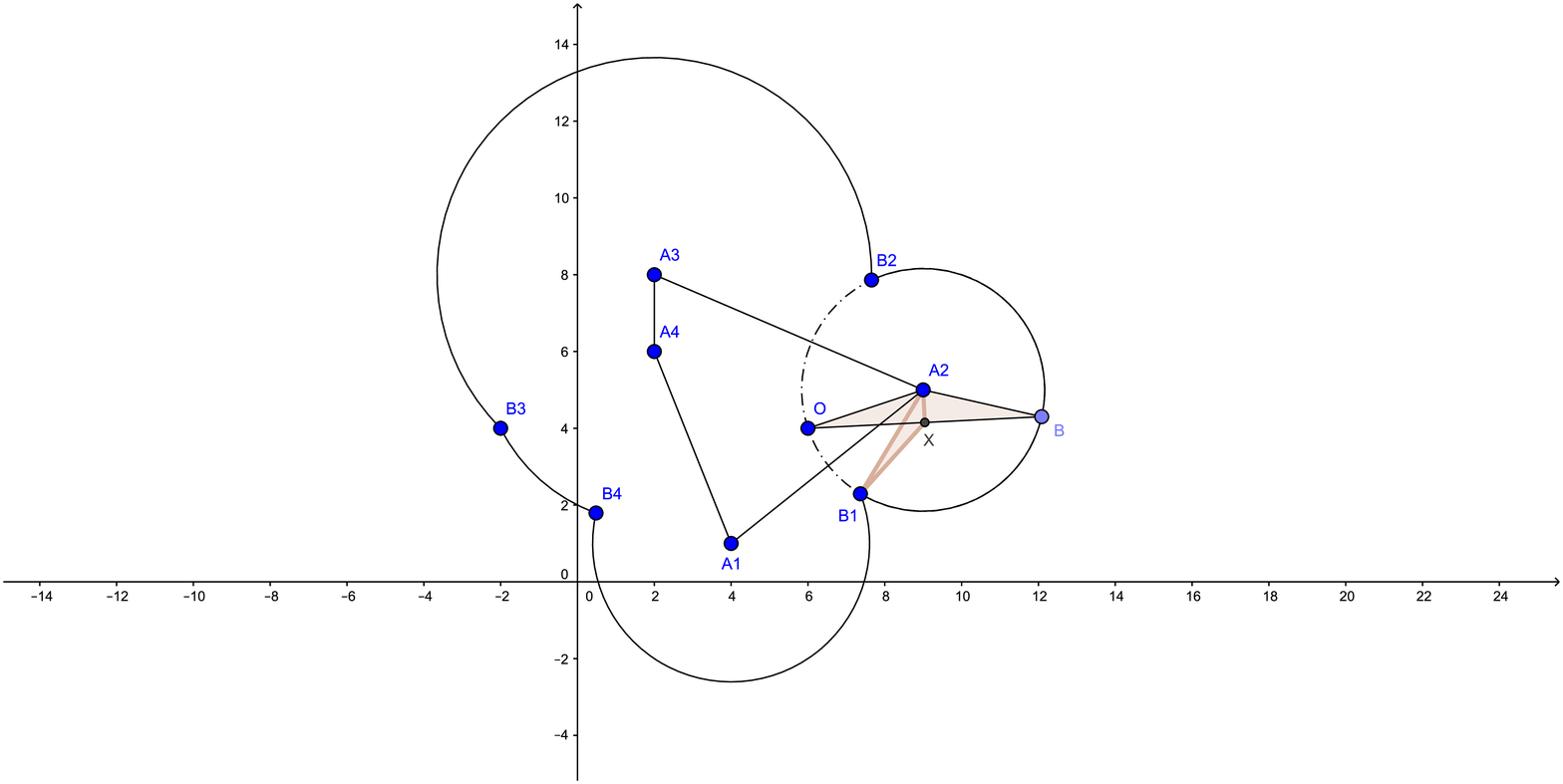}
\caption{The polygon and its focal sets: $K$ and $M$.}
\end{figure}


\end{document}